\documentclass[12pt]{article}
\usepackage{harvard,amsmath,amssymb,latexsym,graphics,epsfig,psfrag,float}

\citationmode{abbr}

\DeclareMathOperator{\corr}{Corr}

\begin{document}


\title{``Pre-conditioning'' for feature selection and  regression  in high-dimensional problems}

\author{
{\sc Debashis Paul}
\thanks{Depts. of Statistics,
Univ. of California, Davis. debashis@wald.ucdavis.edu}\\
{\sc Eric Bair}
\thanks{Depts. of Statistics, Stanford Univ., CA
94305, ebair@stat.stanford.edu}\\
{\sc Trevor Hastie}\thanks{Depts. of Statistics and Health, Research
\&
  Policy,  Stanford
    Univ., CA 94305. hastie@stat.stanford.edu}\\
{\sc Robert Tibshirani}\thanks{Depts. of Health, Research \&
  Policy, and Statistics,
    Stanford Univ, tibs@stat.stanford.edu}
}

\maketitle

\begin{abstract}
 We consider regression problems where the number of predictors greatly
  exceeds the number of observations.
We propose a method for variable selection
that first estimates the regression function, yielding  a
``pre-conditioned'' response variable. The primary method used for
this initial regression is supervised principal components. Then we
apply a standard  procedure such as forward stepwise selection or
the LASSO to the pre-conditioned response variable. In a number of
simulated and real data examples, this two-step procedure
outperforms forward stepwise selection or the usual LASSO (applied
directly to the raw outcome). We also show that under a certain
Gaussian latent variable model, application of the LASSO to the
pre-conditioned response variable is consistent as the number of
predictors and observations increases. Moreover, when the
observational noise is rather large, the suggested procedure can
give a more accurate estimate than LASSO.
We illustrate our method
 on some real problems, including survival analysis with microarray data.

\end{abstract}

\section{Introduction}

In this paper we consider the problem of fitting linear (and other
related) models to data for which the number of features $p$
greatly exceeds  the number of samples $n$. This problem occurs
frequently in genomics, for example in microarray studies in which
$p$ genes are measured on $n$ biological samples.

The problem of model selection for data where number of variables is
typically comparable or much larger than the sample size has
received a lot of attention recently. In particular, various
penalized regression methods are being widely used as means of
selecting the variables having nonzero contribution in a regression
model. Among these tools the $L_1$ penalized regression or LASSO
(\citeasnoun{Ti96}) is one of the  most popular techniques. The
Least Angle Regression (LAR) procedure \citeasnoun{LARS} provides a
method for fast computation of LASSO solution in regression
problems. \citeasnoun{OPT00} derived the optimality conditions
associated with the LASSO solution. \citeasnoun{DE2003} and
\citeasnoun{donoho2004} proved some analytical properties of the
$L_1$ penalization approach for determining the sparsest solution
for an under-determined linear system. Some statistical properties
of the LASSO-based estimator of the regression parameter have been
derived by \citeasnoun{KF2000}. In the context of high-dimensional graphs,
 \citeasnoun{MB2006} showed that the
variable selection method based on lasso  can be consistent if the underlying model
satisfies some conditions.
 Various other model selection criteria have been proposed in
high dimensional regression problems. \citeasnoun{FL2001} and
\citeasnoun{SY2002} gave surveys of some of these methods.

However, when the number of variables ($p$) is much larger than the
number of observations (precisely $p_n \sim c n^\xi$ for some $\xi
\in (0,1)$) \citeasnoun{meinshausen2005} showed that the convergence
rate of risk of the LASSO estimator can be quite slow.
For finite-dimensional problems, \citeasnoun{zou2005} 
 found a necessary condition for the covariance
matrix of the observations, without which  the LASSO variable
selection approach is inconsistent.
\citeasnoun{ZY2006} derived a related result for ther $p >N$ case.

Various modifications to LASSO have been proposed to ensure that on
one hand, the variable selection process is consistent and on the
other, the estimated regression parameter has a fast rate of
convergence. \citeasnoun{FL2001} proposed the Smoothly Clipped
Absolute Deviation (SCAD) penalty for variable selection.
\citeasnoun{FP2004} discussed the asymptotic behavior of this and
other related penalized likelihood procedures when the
dimensionality of the parameter is growing. \citeasnoun{zou2005}
proposed a non-negative Garrote-type penalty (that is re-weighted by
the least squares estimate of the regression parameter) and showed
that this estimator has adaptivity properties when $p$ is fixed.
\citeasnoun{meinshausen2005} proposed a relaxation to the LASSO
penalty after initial model selection to address the problem of high
bias of LASSO estimate when $p$ is very large.

All of these methods try to solve two problems at once: 1)  find a
good predictor $\hat y$ and 2)  find a (hopefully small) subset of
variables to form the basis for this prediction. When $p \gg n$,
these problems are especially difficult. In this paper we suggest
that they should be solved separately, rather than both at once.
Moreover, the method we propose utilizes the correlation structure
of the predictors, unlike most of the methods cited. We propose a
two-stage approach :
\begin{description}
\item  {(a)} find a consistent predictor $\hat y$ of the  true response,
\item{(b)} using the {\em pre-conditioned}  outcome $\hat y$, apply a model
fitting procedure (such as forward stagewise selection or the LASSO)
to the data $(\mathbf{x}, \hat y)$.
\end{description}
In this paper we show that the use of $\hat y$ in place of $y$ in
the model selection step (b) can mitigate the effects of noisy
features on the selection process under the setting of a
\textit{latent variable model} for the response, when the number of
predictor variables that are associated with the response grows at a
slower rate than the number of observations, even though the nominal
dimension of the predictors can grow at a much faster rate.

This paper is organized as follows. In section \ref{sec.precond} we
define the pre-conditioning method and give an example from a latent
variable model. Section \ref{sec.example} discusses a real example
from a kidney cancer microarray study, and application of the idea
to other settings such as survival analysis. In section
\ref{sec.asymptotics} we give details of the latent variable model,
and show that the LASSO applied to the pre-conditioned response
yields a consistent set of predictors, as the number of features and
samples goes to infinity. Finally in section \ref{sec.other} we
discuss and illustrate the pre-conditioning idea for classification
problems.

\section{Pre-conditioning}
\label{sec.precond}

Suppose that the feature measurements are
$x_i=(x_{i1}, x_{i2}, \ldots x_{ip})$ and  outcome values $y_i$, for
$i=1,2,\ldots n$. Our basic model has the form
\begin{eqnarray}
\mathbb{E}(y_i | x_i)= \theta_0+\sum_{j=1}^p x_{ij}\theta_j,
~~~~i=1,2,\ldots,n
\label{eqn:lm}
\end{eqnarray}
Two popular methods for fitting this model are forward stepwise
selection (FS) and the LASSO \citeasnoun{Ti96}.  The first method
successively enters the variable that most reduces the residual sum
of squares, while the second minimizes the penalized criterion
\begin{eqnarray}
J(\theta,\mu)=\sum_i(y_i -\theta_0+\sum_{j=1}^p\theta_j x_{ij})^2
+\mu\sum_{j=1}^p|\theta_j|.
\end{eqnarray}
\citeasnoun{LARS} develop the least angle regression (LAR)
algorithm, for fast computation of the LASSO for all values of the
tuning parameter $\mu \geq 0$.

Usually model selection in the general model (\ref{eqn:lm}) is quite
difficult when $p \gg n$, and our simulations confirm this. To get
better results we may need further assumptions about the underlying
model relating $y_i$ to $x_i$. In this paper, we assume that $y_i$
and $x_i$ are connected via  a low-dimensional {\em latent variable
model}, and use a method that we shall refer to as {\em
pre-conditioning}  to carry out model selection. In this  approach,
we first find a consistent estimate $\hat y_i$ by utilizing the
latent variable structure, and then apply a fitting procedure such
as forward stepwise regression or the LASSO to the data $(x_i, \hat
y_i), i=1,2,\ldots n$. The main technique that we consider for the
initial  pre-conditioning step is supervised principal components
(SPC) (\citeasnoun{BT2004}, \citeasnoun{BHPT2006}). This method  works as follows:
\begin{description}
\item {a)} we select the  features whose individual correlation with
the outcome is large,
\item{b)} using just these features, we compute the principal components
of the matrix of features,  giving  $\hat V_1,
\hat V_2, \ldots \hat V_{\min\{N,p\}}$. The prediction $\hat y_i$ is
the least squares regression of $y_i$ on the first $K$ of these
components.
\end{description}
Typically we use just the first or first few supervised principal
components. \citeasnoun{BHPT2006} show that under an assumption
about the sparsity of the
population principal components, as $p,n \rightarrow \infty$,
supervised principal components gives consistent estimates for the regression
coefficients while the usual principal components regression does
not. We give details of this model in section \ref{sec.asymptotics}, and
provide a simple example next.

\subsection{Example: latent variable model}

The following example shows the main idea n this paper. Consider a model of the form:

\begin{equation}
\label{eq:2}
  Y=\beta_0 + \beta_1 V +\sigma_1 Z
\end{equation}
In addition, we have  measurements on a set of features $X_j$
indexed by $j\in {\cal A}$, for which

\begin{equation}
  X_{j}= \alpha_{0j} + \alpha_{1j} V + \sigma_0 e_{j}, \quad j\in 1,\ldots,p.
\label{eq:latent}
\end{equation}

The  quantity $V$ is an unobserved or {\em latent} variable. The set
${\cal A}$  represents the important features
(meaning that $\alpha_{1j} \neq 0$, for $j \in {\cal A}$) for predicting
$Y_i$. The errors $Z_i$ and $e_{ij}$ are assumed to have mean zero
and are independent of all other random variables in their
respective models. All random variables $(V, Z, e_{j})$ have a
standard Gaussian distribution.

\subsection{Example 1}
\begin{figure}
\begin{psfrags}
\psfrag{corr(u,x)}{corr(V, X)}
\psfrag{corr(y,x)}{corr(Y,X)}
\psfrag{corr(uhat,x)}{corr(Vhat, X)}
\centerline{\epsfig{file=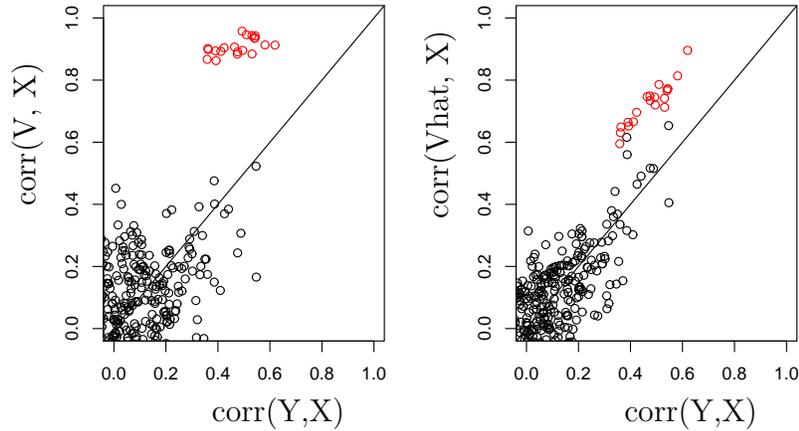,width=.8\textwidth}}
\end{psfrags}
\caption[ex1]{\small\em Results for simulated data. Left panel shows
the correlation between the true latent  variable $V$ and gene
expression $X$ for each of the genes plotted against the correlation
between $Y$ and gene expression. The truly non-null genes are shown
in red. The right panel is the same, except that the estimated
latent  variable $\hat V$ (from supervised principal components) replaces $V$.
We see that  correlation with either the true or estimated latent factor
does a better job  at isolating the  truly non-null genes.}
\label{ex1}
\end{figure}

For illustration, we generated data on $p=500$ features and $n=20$
samples, according to this model, with $\beta_1=2$, $\beta_0=0$,$\alpha_{0j}=0,
\alpha_{1j}=1$,
$\sigma_1 = 2.5$,  ${\cal A} = \{1,2,\ldots 20\}$. Our goal is to predict
$Y$ from $X_{1}, X_{2}, \ldots X_{p}$, and in the process,
discover the fact that only the first 20 features are relevant. This
is a difficult problem. However if we guess (correctly)  that
the data were generated from model (\ref{eq:latent}), our task is made
easier. The left panel of Figure \ref{ex1} shows the correlations
$\corr(V, X_j)$  plotted versus $\corr(Y, X_j)$ for each feature
$j$. The first 20 features are plotted in red, and can be
distinguished much more easily on the basis of $\corr(V,X_j)$ than
$\corr(Y,X_j)$. However this requires knowledge of the underlying
latent factor $V$, which is not observed.

The right panel shows the result when we
instead estimate $V_i$ from the data, using the first supervised
principal component. We see that the correlations of each feature with the
estimated latent factor also distinguishes the relevant from the
irrelevant features.

Not surprisingly, this increased correlation leads to improvements
in the performance of selection methods, as shown in Table
\ref{tab1}. We applied four selection methods to the 20 simulated
data sets from this model:  FS: simple forward stepwise regression;
SPC/FS: forward stepwise regression applied to the pre-conditioned
outcome from supervised principal components; LASSO, and SPC/LASSO:
LASSO applied to pre-conditioned outcome from supervised principal
components.   The table shows the average number of good variables
selected among  the first 1,2,5,10, and 20 variables selected, and
the corresponding  test errors. Pre-conditioning clearly helps 
both forward selection and the lasso.

\begin{table}
\begin{center}
\begin{tabular}{l|rrrr|rrrr}
Method & \multicolumn{4}{c}{Mean \# of good variables,}\\
 & \multicolumn{4}{c}{ when selecting first: } &
\multicolumn{4}{c}{Test error when selecting first:}\\
 &1  &5  & 10  & 20  & 1   &5  & 10 & 20 \\
\hline
 FS         &  0.82&  0.98&  1.12  & 1.58 &   267.36& 335.4 & 353.52& 357.07\\
 SPC/FS  &  0.94&  2.66 & 2.86 &  3.12  &  241.88& 229.47& 231.52& 232.28\\
 LASSO      &  0.88&  2.05  & 3.17 &  3.29 &   206.54 & 184.56
 & 186.71& 205.85\\
 SPC/LASSO& 0.92&  4.21 & 7.75 &  9.71 &   212.23& 197.07& 183.04& 178.19\\
\end{tabular}
\end{center}
\caption{\em Four selection methods to the 20 simulated data sets
from the model of Example 1.
Shown are  the number of good variables selected among  the first
1,2,10, and 20 variables selected, and the corresponding  test
errors. Pre-conditioning clearly helps in both cases,
and the lasso outperforms forward selection.}
 \label{tab1}
\end{table}

\subsection{Example 2.}
The second example was suggested by a referee. 
 It is somewhat artifical
but exposes an important
assumption that is made by our procedure.
We define random variables
$(Y, X_1, X_2, X_3)$  having a Gaussian distribution with mean zero
and inverse  covariance matrix
$$
\Sigma^{-1}= \left( \begin{array}{cccc}
2 & 1 & 1 &1 \\
1 & 2 & 0 &1 \\
1 & 0 & 2 &1 \\
1 & 1 & 1 & 2 \end{array} \right).
$$
We define 297 additional predictors that are $N(0,1)$.
The  population regression coefficient is
$\beta=(-1, -1, -1, 0,0,\ldots)$ while  the (marginal)
correlation of each predictor with $Y$ is
$\rho=(-0.5, -0.5,0,  0,0,\ldots)$. Hence $X_3$ has zero marginal
correlation with $Y$ but has a non-zero partial correlation with $Y$,
(since $(\Sigma^{-1})_{14}=1$).
The number of  good variables when selecting the first 1,2,3 or 4 predictors
is shown in Table \ref{tab:sim2}.

\begin{table}
\begin{center}
\begin{tabular}{l|rrrr}|
Method & \multicolumn{4}{c}{Mean \# of good variables.} \\
& \multicolumn{4}{c}{ when selecting first: } \\
 &1  &2  & 3  & 4 \\
\hline
 LASSO      &  1.0&  2.0  & 3.0 &  3.0 \\
 SPC/LASSO& 1.0&  2.0 & 2.0 &  2.0 \\
\end{tabular}
\end{center}
\caption{\em Performance of LASSO and  pre-conditioned LASSO in the
second simulation example.} \label{tab:sim2}
\end{table}
We see that the LASSO enters the 3 good predictors first in every simulation,
while the pre-conditioned version ignores the 3rd predictor.
Supervised principal components screens out this  predictor,
because it is marginally independent of $Y$.

Pre-conditioning  with supervised principal components assumes that
any important predictor (in the sense of having significantly large
nonzero regression coefficient) will also have a substantial
marginal correlation with the outcome.  This need not be true in
practice, but we believe it will often be a good working hypothesis
in many practical problems.

\subsection{Example 3.}
Our third simulation study compares  the lasso to the pre-conditioned lasso,
in a more neutral setting. We generated 1000 predictors,  each  having a $N(0,1)$
distribution marginally. The first 40 predictors had a pairwise correlation of 0.5, while the
remainder were uncorrelated.

The outcome was generated as
\begin{eqnarray}
Y= \sum_{j=1}^{40} \beta_j X_j +\sigma Z
\end{eqnarray}
with $Z, \beta_j \sim N(0,1)$ and $\sigma=5$.
Hence the outcome is only a function of the first 40 (``good'')  predictors.

We generated 100 datasets from this model:
the average number of good variables selected by the lasso and pre-conditioned
lasso is shown in Table \ref{tab:sim3}.
\begin{table}
\begin{center}
\begin{tabular}{l|rrrr}|
Method & \multicolumn{4}{c}{Mean \# of good variables.} \\
& \multicolumn{4}{c}{ when selecting first: } \\
 &5  &10  & 20  & 50 \\
\hline
LASSO &     2.92&  5.88 & 9.04 & 9.16\\
SPC/LASSO &     2.49&  5.13& 10.32& 19.73\\
\end{tabular}
\end{center}
\caption{\em Performance of LASSO and  pre-conditioned LASSO in the
third simulation example.} \label{tab:sim3}
\end{table}
Note that with just $n=50$ samples, the maximum number of predictors in the
model is also 50.
While neither method  is  successful at isolating the bulk of the  40 good predictors,
the pre-conditioned lasso finds twice as many  good predictors as the lasso
in the  full model.

\section{Examples}\label{sec.example}

\subsection{Kidney cancer data}

\citeasnoun{Zhaoetal2005} collected gene expression  data on $14,814$
genes from 177 kidney patients. Survival times (possibly censored)
were also measured for each patient,
as well as a number of clinical predictors including the grade of the
tumor: 1 (good)  to 4 (poor).

\begin{figure}
\centerline{\epsfig{file=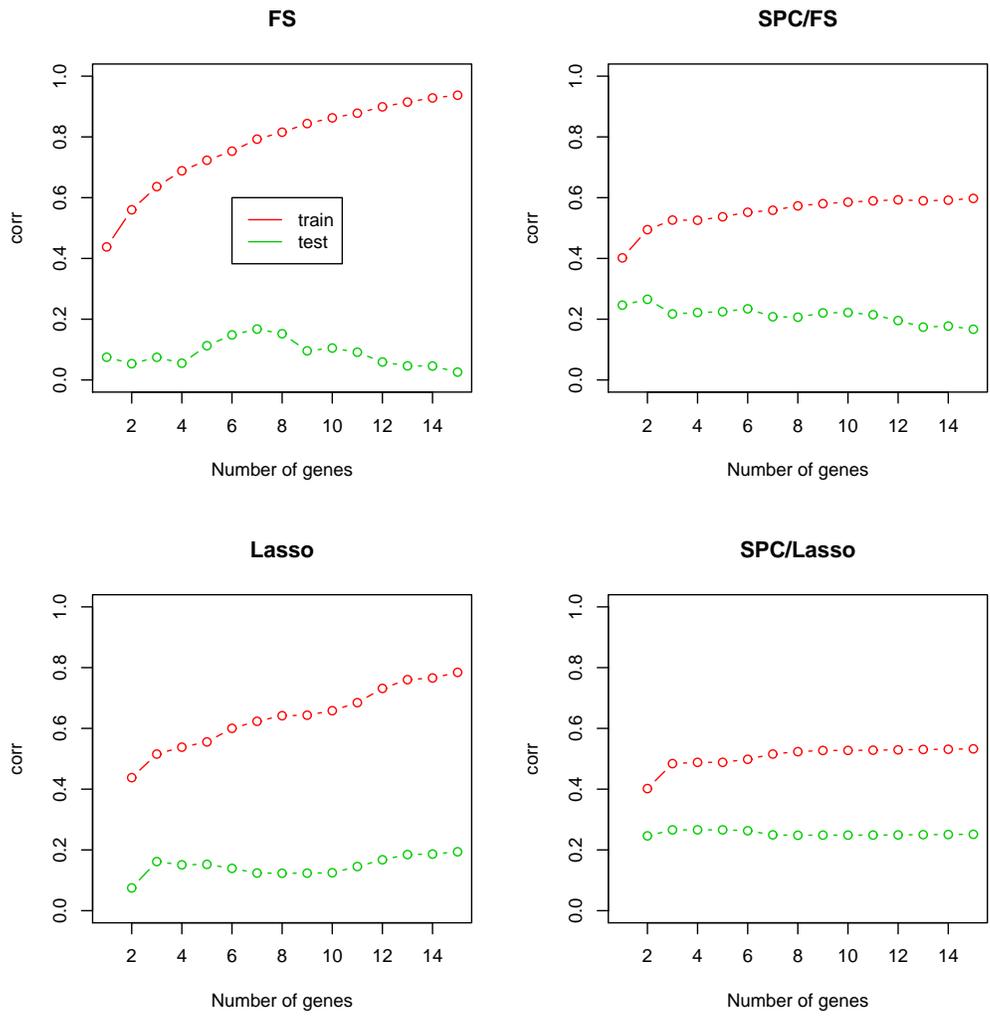,width=\textwidth}}
\caption[brooks.grade]{\small\em Kidney cancer data: predicting
tumor grade. Correlation of different predictors with the true
outcome, in the training and test sets, as more and more  genes are
entered. } \label{brooks.grade}
\end{figure}

The data were split into 88 samples to form the  training set and
the remaining 89 formed the test set. For illustration, in this
section we try to predict grade from gene expression. In the next
section we predict survival time (the primary outcome of interest)
from gene expression. Figure \ref{brooks.grade} shows the training
and test set correlations between grade  and its prediction from
different methods. We see that for both forward selection and the
LASSO, use of the supervised principal component prediction $\hat y$
as the outcome variable (instead of $y$ itself) makes the procedure
less greedy in the training set and yields higher correlations in
the test set. While the correlations in the test set are not
spectacularly high, for SPC/FS and SPC/LASSO they do result in a
better predictions in the test set.

\subsection{Application to other regression settings}

Extension of our proposal to other kinds of regression outcomes  is very
simple. The only change is in step (a) of
supervised principal components algorithm, where we replace the correlation
by an appropriate measure of association. In particular, the likelihood
score statistic is an attractive choice.

\subsection{Survival analysis}

Perhaps the most common version of the $p > n$ regression problem in
genomic studies is survival analysis, where the outcome is patient
survival (possibly censored). Then  we use  the partial likelihood
score statistic from Cox's proportional hazards score statistic (see
Chapter 4 of \citeasnoun{KP80}), in step (a) of supervised principal
components. After that, we can (conveniently) use the usual least
squares version of FS or LASSO in step (2) of the modeling process.
Hence the computational advantages of the least angle regression
algorithm can be exploited.

\begin{figure}
\begin{psfrags}
\psfrag{usual}{FS} \psfrag{SPC}{SPC/FS}
\centerline{\epsfig{file=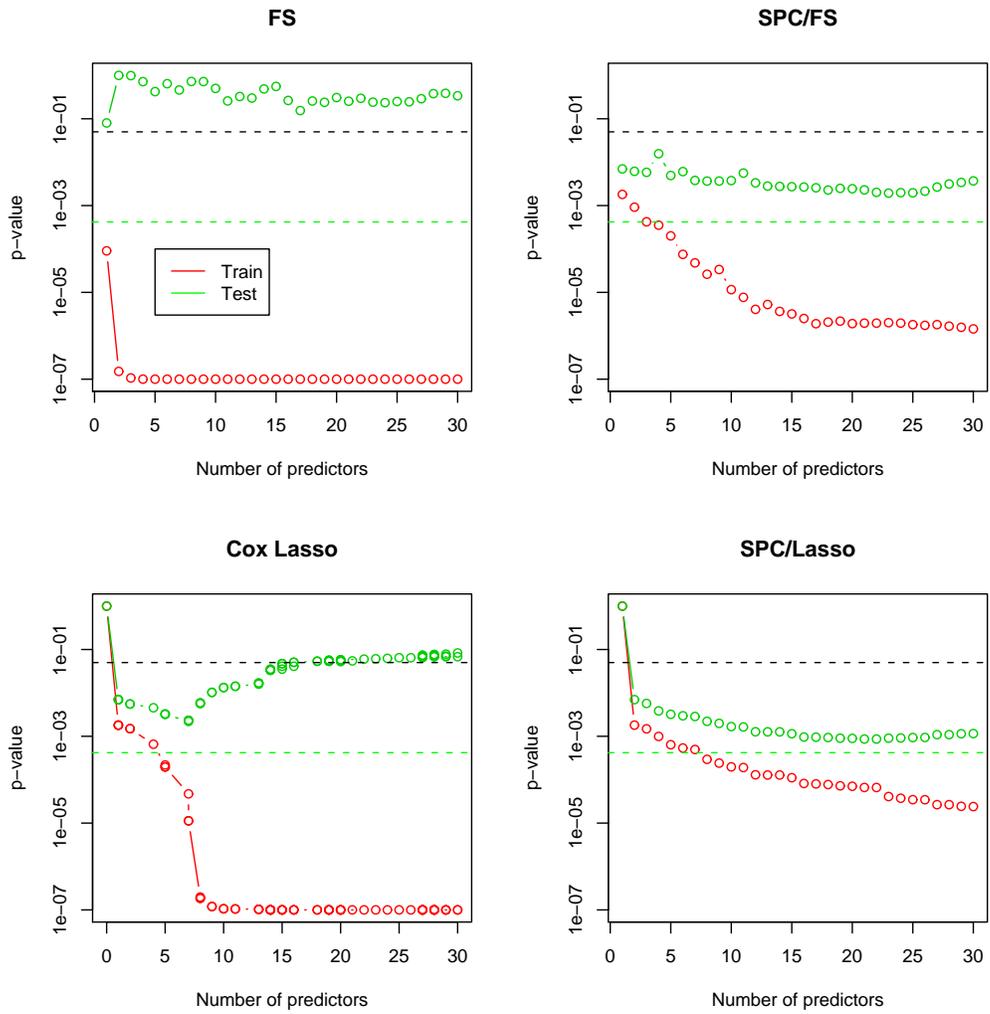,width=\textwidth}}
\end{psfrags}
\caption[figbrooks]{\small\em Kidney cancer  data:  predicting
  survival time. Training set p-values (red) and test set p-values
  (green) for four different
 selection methods as more and more  genes are entered.
Horizontal broken  lines are drawn at 0.05 (black) and  the test set p-value
for the supervised principal component predictor 0.00042 (green).}
\label{figbrooks}
\end{figure}

Figure \ref{figbrooks} shows the result of applying forward stepwise
Cox regression (top left panel), forward stepwise selection applied to the
SPC predictor (top right panel), LASSO for the Cox model (bottom
left panel) and LASSO applied  to the SPC
predictor (bottom right panel).
The bottom left panel was computed using the {\tt glmpath} R package
of \citeasnoun{PH2006}, available in the CRAN collection.
  In each case we obtain a predictor
$\hat y$, and then use $\hat y$ as a covariate in a Cox model, in
either the training or test set. The resulting p-values from these Cox
models are shown in the figure. We see that  forward stepwise Cox
regression tends to overfit in the training set, and hence the
resulting test-set p-values are not significant. The two stage
SPC/FS procedure fits more slowly in the training set, and hence
achieves smaller p-values in the test set. 
``SPC/LASSO'' , the LASSO applied to the pre-conditioned
response from  supervised principal components, performs best
and is  also computationally convenient:
it uses the fast LAR algorithm for the lasso, applied to the
pre-conditioned response variable.

The horizontal green line shows the test set p-value of the
supervised principal component predictor. We see that  the first 10
or 15 genes chosen by the LASSO have captured the  signal in this
predictor.

We have used the pre-conditioning procedure in real microarray
studies. We have found that it is useful to report to investigators
not just the best 10 or 15 gene model, but also any genes that have
high correlation with this set. The enlarged set can be useful in
understanding the underlying biology in experiment, and also for
building assays  for future clinical use. A given gene might not be
well measured on a microarray for a  variety of reasons, and hence
it is useful to identify surrogate genes that may be used in its
place.

\begin{figure}
\centerline{\epsfig{file=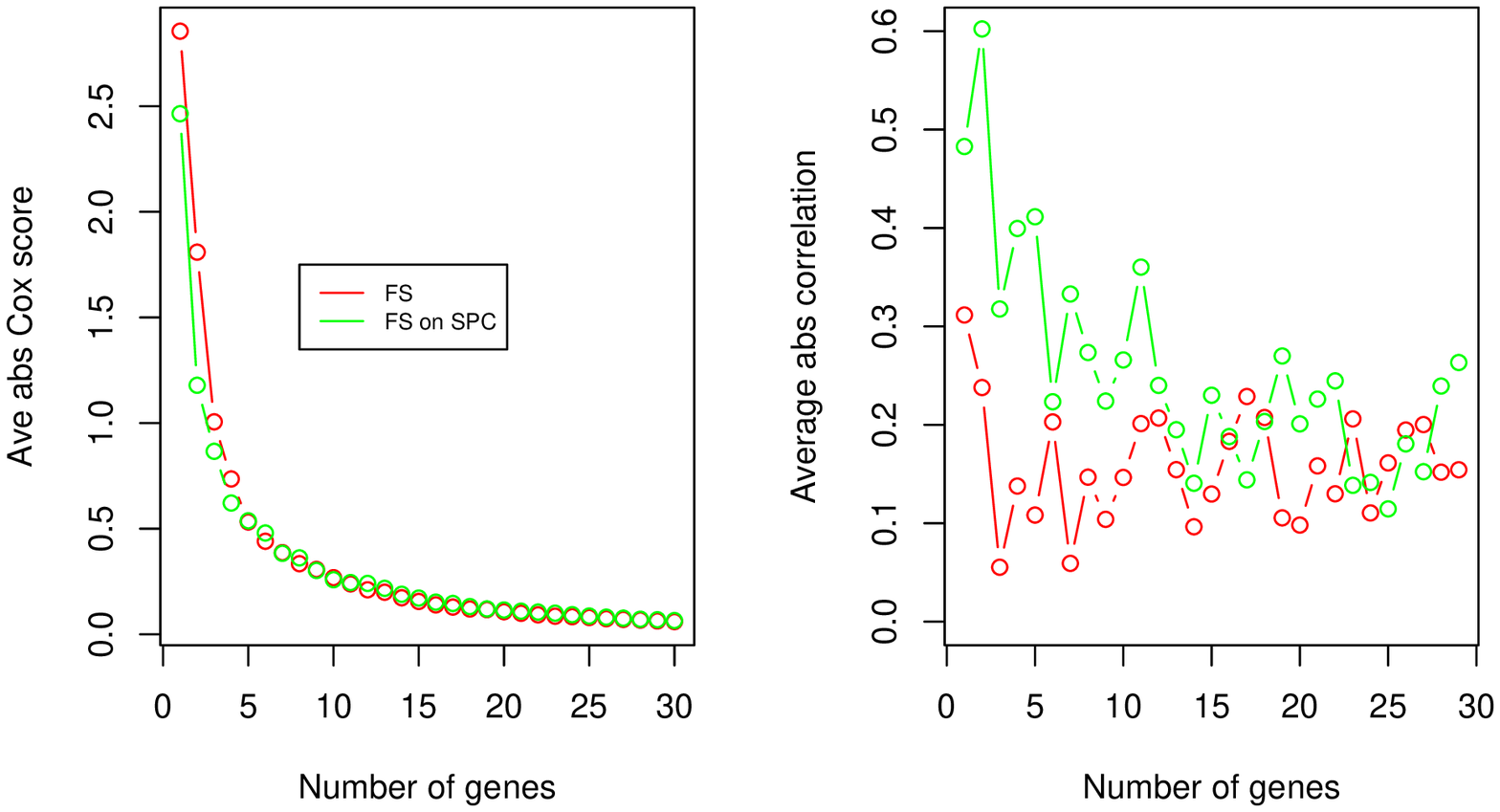,width=.8\textwidth}}
\caption[figbrooks3]{\small\em Kidney cancer data:  predicting
survival time. Left panel shows the average absolute Cox score
of the first $k$  genes entered
by  forward stepwise selection (red) and the pre-conditioned version (green),
as $k$ runs from 1 to 30. The right panel
shows the average absolute pairwise correlation of the genes for both methods.}
\label{figbrooks3}
\end{figure}

Figure \ref{figbrooks3}  shows the  average absolute Cox score
of the first $k$  features entered by  forward stepwise selection
(red) and the pre-conditioned version (green), as $k$ runs from 1 to
30. The right panel shows the average absolute pairwise correlation
of the genes for both methods. We see that the methods enter
features of about the same strength, but pre-conditioning enters
genes that are more highly correlated with one another.

\section{Asymptotic analysis}\label{sec.asymptotics}

In this section we lay down a mathematical formulation of the
problem and pre-conditioning procedure in the context of a latent
factor model for the response. We show that the procedure combining
SPC with LASSO, under some assumptions about the correlation
structure among the variables, leads to asymptotically consistent
variable selection in the Gaussian linear model setting. We consider
the class of problems where one observes $n$ independent samples
$(y_i,\mathbf{x}_i)$ where $y_i$ is a one dimensional response and
$\mathbf{x}_i$ is a $p$-dimensional predictor. Individual
coordinates of the vector $\mathbf{x}_i$ are denoted by $x_{ij}$
where the index $j \in \{1,\ldots,p\}$ correspond to the $j$-th
predictor. We denote the $n \times p$ matrix $((x_{ij}))_{1\leq i
\leq n,1\leq j \leq p}$ by $\mathbf{X}$ and the vector
$(y_i)_{i=1}^n$ by $Y$. Henceforth, unless otherwise stated, we do
not make a distinction between the realized value $(Y,\mathbf{X})$
and the random elements (namely, the response and the $p$
predictors) that they represent.

The interest is in identifying the set of predictors $X_j$ which are
(linearly) related to $Y$. A regression model will be of the form
$\mathbb{E}(Y|\mathbf{x}) = \theta^T \mathbf{x}$ for some $\theta
\in \mathbb{R}^p$. Here we assume that the joint distribution of
$\mathbf{X}$ is Gaussian with zero mean and covariance matrix
$\Sigma \equiv \Sigma_p$. The relationship between $Y$ and
$\mathbf{X}$ is assumed to be specified by a latent component model
to be described below.

\subsection{Model for $\mathbf{X}$}

Suppose that the spectral decomposition of $\Sigma$ is given by
$\Sigma = \sum_{k=1}^p \ell_k \mathbf{u}_k \mathbf{u}_k^T$, where
$\ell_1 \geq \ldots \geq \ell_p \geq 0$ and $\mathbf{u}_1,\ldots,
\mathbf{u}_p$ form an orthonormal basis of $\mathbb{R}^p$. We
consider the following model for $\Sigma$.

Assume that there exists an $M \geq 1$ such that
\begin{equation}\label{eq:signal_plus_noise}
\ell_k = \lambda_k + \sigma_0^2, ~~k=1,\ldots,M, ~\mbox{and}~ \ell_k = \sigma_0^2, ~~
k=M+1,\ldots,p,
\end{equation}
where $\lambda_1 \geq \ldots \geq \lambda_M > 0$ and $\sigma_0 > 0$.
This model will be referred to as the ``noisy factor model''. To see
this, notice that under the Gaussian assumption the matrix
$\mathbf{X}$ can be expressed as
\begin{equation}\label{eq:NFM_model}
\mathbf{X} = \sum_{k=1}^M \sqrt{\lambda_k} \mathbf{v}_k
\mathbf{u}_k^T + \sigma_0 \mathbf{E}
\end{equation}
where $\mathbf{v}_1,\ldots,\mathbf{v}_M$ are i.i.d. $N_n(0,I)$
vectors (the factors), and $\mathbf{E}$ is an $n \times p$ matrix
with i.i.d. $N(0,1)$ entries, and is independent of
$\mathbf{v}_1,\ldots,\mathbf{v}_M$. This matrix is viewed as a noise
matrix.

In the analysis presented in this paper throughout we use
(\ref{eq:NFM_model}) as the model for $\mathbf{X}$, even though it
can be shown that the analysis applies even in the case where
$\ell_{K+1},\ldots,\ell_p$ are decreasing and sufficiently well
separated from $\ell_1,\ldots,\ell_K$.


\subsection{Model for $Y$}

Assume the following regression model for $Y$. Note that this is a
more general version of (\ref{eq:2}), even though we assume that $Y$
has (unconditional) mean 0.
\begin{equation}\label{eq:Y_model}
Y = \sum_{k=1}^K \beta_k \mathbf{v}_k + \sigma_1 Z,
\end{equation}
where $\sigma_1 > 0$, $1\leq K \leq M$, and $Z$ has $N_n(0,I)$
distribution and is independent of $\mathbf{X}$.


\subsection{Least squares and feature selection}

We derive expressions for the marginal correlations between $Y$ and
$X_j$, for $j=1,\ldots,p$ and the (population) least squares
solution, viz. $\theta := \arg\min_{\zeta}\mathbb{E}
\parallel Y - \mathbf{X}\zeta\parallel_2^2$,
in terms of the model parameters. Let ${\cal P} := \{1,\ldots,p\}$.
The marginal correlation between $\mathbf{X} = (X_j)_{j=1}^p$ and
$Y$ is given by
\begin{equation}\label{eq:marginal}
\Sigma_{{\cal P} y} := (\mathbb{E}(X_j Y))_{j=1}^p = \sum_{k=1}^K
\beta_k \sqrt{\lambda_k} \mathbf{u}_k.
\end{equation}
The population regression coefficient of $Y$ on $\mathbf{X}$, is
given by
\begin{eqnarray}\label{eq:regression}
\theta = \Sigma^{-1} \Sigma_{{\cal P} y} &=& [\sum_{k=1}^M \lambda_k
\mathbf{u}_k \mathbf{u}_k^T + \sigma_0^2 I]^{-1}
[\sum_{k=1}^K \beta_k \sqrt{\lambda_k} \mathbf{u}_k] \nonumber \\
&=& \left[\sum_{k=1}^M \frac{1}{\lambda_k + \sigma_0^2} \mathbf{u}_k
\mathbf{u}_k^T + \frac{1}{\sigma_0^2} (I - \sum_{k=1}^M \mathbf{u}_k
\mathbf{u}_k^T)\right] [\sum_{k=1}^K \beta_k \sqrt{\lambda_k}
\mathbf{u}_k]  \nonumber\\
&=& \sum_{k=1}^K \beta_k \frac{\sqrt{\lambda_k}}{\lambda_k +
\sigma_0^2} \mathbf{u}_k ~=~ \sum_{k=1}^K \beta_k \ell_k^{-1}
\sqrt{\lambda_k}\mathbf{u}_k.
\end{eqnarray}
Now, define $\mathbf{w}_j = (\sqrt{\lambda_1}
u_{j1},\ldots,\sqrt{\lambda_K} u_{jK})^T$. Let ${\cal D} = \{j :
\parallel \mathbf{w}_j \parallel_2 \neq 0\}$. Observe that
$\Sigma_{jy} = \beta^T \mathbf{w}_j$, and $\theta_j = \beta^T
D_K^{-1} \mathbf{w}_j$, where $D_K =$ diag$(\ell_1,\ldots,\ell_K)$.
So if we define ${\cal B} := \{j : \Sigma_{jy} \neq 0\}$, and ${\cal
A} = \{j : \theta_j \neq 0\}$, then ${\cal B}\subset {\cal D}$ and
${\cal A} \subset {\cal D}$.

This gives rise to the regression model:
\begin{equation}\label{eq:NFM_regression}
Y = \mathbf{X}\theta + \sigma_\varepsilon \varepsilon,
\end{equation}
where
\begin{equation}\label{eq:sigma_error}
\sigma_\varepsilon^2 = \sigma_{yy} - \Sigma_{y{\cal P}}\Sigma^{-1}
\Sigma_{{\cal P}y} = \sigma_1^2 + \sum_{k=1}^K \beta_k^2 -
\sum_{k=1}^K \beta_k^2 \frac{\lambda_k}{\lambda_k + \sigma_0^2} =
\sigma_1^2 + \sigma_0^2 \beta^T D_K^{-1} \beta,
\end{equation}
and $\varepsilon$ has i.i.d. $N(0,1)$ entries and is independent of
$\mathbf{X}$.

Note also that, the population partial covariance between $Y$ and
$\mathbf{X}_C$ given $\mathbf{X}_{\cal D}$ (given by $\Sigma_{y
C|{\cal D}} := \Sigma_{yC} - \Sigma_{y {\cal D}} \Sigma_{{\cal
D}{\cal D}}^{-1} \Sigma_{{\cal D}C}$), for any subset $C \subset
{\cal D}^c$, where ${\cal D}^c := {\cal P} \setminus {\cal D}$, is
0. However the corresponding statement is not true in general if one
replaces ${\cal D}$ by either ${\cal A}$ or ${\cal B}$. Therefore,
ideally, one would like to identify ${\cal D}$. However, it may not
be possible to accomplish this in general when the dimension $p$
grows with the sample size $n$. Rather, we define the feature
selection problem as the problem of identifying ${\cal A}$, while
the estimation problem is to obtain an estimate of $\theta$ from
model (\ref{eq:NFM_regression}).

Observe that, if either $K=1$ or $\lambda_1 = \cdots = \lambda_K$,
then ${\cal A} = {\cal B}$. In the former case we actually have
${\cal A} = {\cal B} = {\cal D}$. In these special cases, the
feature selection problem reduces to finding the set ${\cal B}$,
which may be done (under suitable identifiability conditions) just
by computing the sample marginal correlations between the response
and the predictors and selecting those variables (coordinates) for
which the marginal correlation exceeds an appropriate threshold. The
major assumptions that we shall make here for solving the problem
are that (i) ${\cal A} \subset {\cal B}$, (ii) ${\cal B}$ can be
identified from the data (at least asymptotically), (iii)
cardinality of ${\cal B}$ (and hence that of ${\cal A}$) is small
compared to $n$, and (iv) the contribution of the coordinates ${\cal
B}^c$ in the vectors $\mathbf{u}_1,\ldots,\mathbf{u}_K$ is
asymptotically negligible in an $L^2$ sense. If these conditions are
satisfied, then it will allow for the identification of ${\cal A}$,
even as dimension increases with the sample size. We make these (and
other) conditions more precise in Section \ref{sec:assumptions}.

\subsection{SPC as a preconditioner}

The formulation in the previous section indicates that one may use
some penalized regression methods to estimate the regression
parameter $\theta$ from the model (\ref{eq:NFM_regression}).
However, standard methods like LASSO do not use the covariance
structure of the data. Therefore if one uses the underlying
structure for $\Sigma$, and has good estimates of the parameters
$(\mathbf{u}_k, \ell_k)$, then one can hope to be able to obtain a
better estimate $\theta$, as well as identify ${\cal A}$ as $n \to
\infty$.

We focus on (\ref{eq:NFM_model}) and (\ref{eq:Y_model}).
In general it is not possible to eliminate the contribution of
$\mathbf{E}$ entirely from an estimate of $\mathbf{v}_k$, even if we
had perfect knowledge of $(\mathbf{u}_k,\ell_k)$. To understand
this, note that, the conditional distribution of $\mathbf{v}_k$
given $\mathbf{X}$ is the same as the conditional distribution of
$\mathbf{v}_k$ given $\mathbf{X}\mathbf{u}_k$. The latter
distribution is normal with mean $\frac{\sqrt{\lambda_k}}{\ell_k}
\mathbf{X}\mathbf{u}_k$ and covariance matrix
$\frac{\sigma_0^2}{\ell_k} I_n$. This means that any reasonable
procedure that estimates the parameters $(\mathbf{u}_k,\ell_k)$ can
only hope to reduce the effect of the measurement noise in $Y$, viz.
$\sigma_1 Z$.

Keeping these considerations in mind, we employ a two stage
procedure described in the following section for estimating
$\theta$. In order to fit the model (\ref{eq:NFM_regression}) using
SPC procedure, it is necessary to estimate the eigenvectors
$\mathbf{u}_k$, $k=1,\ldots,M$. When $\frac{p}{n}$ is large (in the
sense that the fraction does not converge to 0 as $n \to \infty$),
in general it is not possible to estimate $\mathbf{u}_k$
consistently. However, if $\mathbf{u}_k$ are sparse, in the sense of
having say $q$ non-zero components, where $\frac{q}{n} \to 0$, then
\citeasnoun{BHPT2006} showed that under suitable identifiability
conditions, it is possible to get asymptotically consistent
estimators of $\mathbf{u}_1,\ldots,\mathbf{u}_K$, where the
consistency is measured in terms of convergence of the $L^2$
distance between the parameter and its estimator.

\subsection{Algorithm}

In this section we present the algorithm in detail.

\begin{itemize}
\item[{\bf Step 1}]
Estimate $(\mathbf{u}_1,\ell_1),\ldots,(\mathbf{u}_K,\ell_K)$ by SPC
procedure in which only those predictors $X_j$ whose empirical
correlation with response $Y$ is above a threshold $\tau_n$ are used
in the eigen-analysis. Call these estimates
$\{\widetilde{\mathbf{u}}_k, \widetilde \ell_k\}_{k=1}^K$.


\item[{\bf Step 2}]
Let $\widetilde P_K := ~Proj~(\widehat V_1,\ldots,\widehat V_K)$ be
the projection onto $\widehat V_1,\ldots,\widehat V_K$, where
$\widehat V_k := \frac{1}{\sqrt{\widetilde\ell_k}}
\mathbf{X}\widetilde{\mathbf{u}}_k$ is the $k$-th principal
component of the predictors (under the SPC procedure). Define
$\widetilde Y = \widetilde P_K Y$.

\item[{\bf Step 3}]
Estimate $\theta$ from the linear model $\widetilde Y = \mathbf{X}
\theta +$ error, using the LASSO approach with penalty $\mu_n
>0$.
\end{itemize}

Since by definition $\frac{1}{n}
\langle\mathbf{X}\widetilde{\mathbf{u}}_k,
\mathbf{X}\widetilde{\mathbf{u}}_{k'} \rangle = \widetilde \ell_k
\delta_{kk'}$, it follows that
\begin{equation}\label{eq:P_tilde_K}
\widetilde P_K =
~Proj~(\mathbf{X}\widetilde{\mathbf{u}}_1,\ldots,\mathbf{X}\widetilde{\mathbf{u}}_K)
= \sum_{k=1}^K
\frac{1}{\parallel\mathbf{X}\widetilde{\mathbf{u}}_k\parallel^2}
(\mathbf{X}\widetilde{\mathbf{u}}_k)
(\mathbf{X}\widetilde{\mathbf{u}}_k)^T = \sum_{k=1}^K
\frac{1}{\widetilde \ell_k}
\frac{1}{n}(\mathbf{X}\widetilde{\mathbf{u}}_k)
(\mathbf{X}\widetilde{\mathbf{u}}_k)^T.
\end{equation}

\subsection{Analysis of the projection}

We present an expansion of the projected response $\widetilde Y :=
\widetilde P_K Y$ that will be useful for all the asymptotic
analyses that follow. Using the representation of $\widetilde P_K$
in (\ref{eq:P_tilde_K}) and invoking (\ref{eq:NFM_model}) and
(\ref{eq:Y_model}), we get
\begin{eqnarray}\label{eq:Y_tilde_bar_approx}
\widetilde Y &=& \sum_{k=1}^K \frac{\beta_k}{\widetilde \ell_k}
\frac{1}{n}\langle \mathbf{X}\widetilde{\mathbf{u}}_k, \mathbf{v}_k
\rangle \mathbf{X}\widetilde{\mathbf{u}}_k +
\sum_{k=1}^K\sum_{k'\neq k}^K \frac{\beta_{k'}}{\widetilde \ell_k}
\frac{1}{n}\langle \mathbf{X}\widetilde{\mathbf{u}}_k,
\mathbf{v}_{k'} \rangle \mathbf{X}\widetilde{\mathbf{u}}_k +
\sigma_1 \sum_{k=1}^K \frac{1}{\widetilde\ell_k} \frac{1}{n} \langle
\mathbf{X} \widetilde{\mathbf{u}}_k, Z \rangle \mathbf{X}
\widetilde{\mathbf{u}}_k \nonumber\\
&=& \sum_{k=1}^K \frac{\beta_k\sqrt{\lambda_k}}{\widetilde \ell_k}
\frac{1}{n}\parallel \mathbf{v}_k \parallel^2 \langle \mathbf{u}_k,
\widetilde{\mathbf{u}}_k \rangle \mathbf{X}\widetilde{\mathbf{u}}_k
+ \sum_{k=1}^K \sum_{l \neq k}^M
\frac{\beta_k\sqrt{\lambda_l}}{\widetilde \ell_k} \frac{1}{n}\langle
\mathbf{v}_l,\mathbf{v}_k \rangle \langle \mathbf{u}_l,
\widetilde{\mathbf{u}}_k \rangle
\mathbf{X}\widetilde{\mathbf{u}}_k \nonumber\\
&& + \sum_{k=1}^K \sum_{k'\neq k}^K \sum_{l=1}^M
\frac{\beta_{k'}\sqrt{\lambda_l}}{\widetilde \ell_k}
\frac{1}{n}\langle \mathbf{v}_l, \mathbf{v}_{k'} \rangle \langle
\mathbf{u}_l, \widetilde{\mathbf{u}}_k \rangle
\mathbf{X}\widetilde{\mathbf{u}}_k \nonumber\\
&& + \sigma_0 \sum_{k=1}^K \sum_{k'=1}^K\frac{\beta_{k'}}{\widetilde
\ell_k} \frac{1}{n} \langle \mathbf{E}\widetilde{\mathbf{u}}_k,
\mathbf{v}_{k'} \rangle \mathbf{X} \widetilde{\mathbf{u}}_k +
\sigma_1 \sum_{k=1}^K \frac{1}{\widetilde \ell_k} \frac{1}{n}
\langle \mathbf{X} \widetilde{\mathbf{u}}_k, Z \rangle \mathbf{X}
\widetilde{\mathbf{u}}_k \nonumber\\
&=& \mathbf{X} \theta + \mathbf{X} \sum_{k=1}^K
\beta_k\sqrt{\lambda_k}(\frac{1}{\widetilde\ell_k}\frac{\parallel
\mathbf{v}_k\parallel^2}{n} \langle
\mathbf{u}_k,\widetilde{\mathbf{u}}_k\rangle
\widetilde{\mathbf{u}}_k - \frac{1}{\ell_k} \mathbf{u}_k) +
\sum_{k=1}^K \sum_{k'\neq k}^K
\frac{\beta_{k'}\sqrt{\lambda_{k'}}}{\widetilde \ell_k}
\frac{\parallel \mathbf{v}_{k'} \parallel^2}{n} \langle
\mathbf{u}_{k'}, \widetilde{\mathbf{u}}_k \rangle
\mathbf{X}\widetilde{\mathbf{u}}_k \nonumber\\
&&  + \sum_{k=1}^K \sum_{k'\neq k}^K
\frac{\beta_{k'}\sqrt{\lambda_k}}{\widetilde \ell_k}
\frac{1}{n}\langle \mathbf{v}_k, \mathbf{v}_{k'} \rangle \langle
\mathbf{u}_k, \widetilde{\mathbf{u}}_k \rangle
\mathbf{X}\widetilde{\mathbf{u}}_k \nonumber\\
&& + \sigma_0 \sum_{k=1}^K \sum_{k'=1}^K
\frac{\beta_{k'}}{\widetilde \ell_k} \frac{1}{n} \langle \mathbf{E}
\widetilde{\mathbf{u}}_k, \mathbf{v}_{k'} \rangle \mathbf{X}
\widetilde{\mathbf{u}}_k + \sigma_1 \sum_{k=1}^K \frac{1}{\widetilde
\ell_k} \frac{1}{n} \langle \mathbf{X} \widetilde{\mathbf{u}}_k, Z
\rangle \mathbf{X} \widetilde{\mathbf{u}}_k + R_n,
\end{eqnarray}
for some vector $R_n \in \mathbb{R}^n$. This is an asymptotically
unbiased regression model for estimating $\theta$ provided
$(\widetilde{\mathbf{u}}_k,\widetilde\ell_k)_{k=1}^K$ is an
asymptotically consistent estimator for
$(\mathbf{u}_k,\ell_k)_{k=1}^K$.



\subsection{Assumptions}\label{sec:assumptions}

In this section we give sufficient conditions for the consistency of
the variable selection aspect of the SPC preconditioning procedure.
The methods of \citeasnoun{zou2005} and \citeasnoun{KF2000} are not
applicable in our situation since the dimension is growing with the
sample size. For most parts, we make assumptions similar to those in
\citeasnoun{MB2006} for the relationship among the variables.

\begin{itemize}
\item[{\bf A1}] The eigenvalues $\lambda_1,\ldots,\lambda_M$ satisfy
\begin{itemize}
\item[(i)]
$\lambda_1 > \ldots > \lambda_K > \lambda_{K+1} \geq \ldots\geq
\lambda_M \geq 0$.
\item[(ii)] $\min_{1\leq k \leq K} (\lambda_k -
\lambda_{k+1})\geq C_0$ for some $C_0 > 0$ (fixed).
\item[(iii)] $\lambda_1 \leq \Lambda_{\max}$ for some
$\Lambda_{\max}$ fixed. Also, $\sigma_0$ is fixed.
\end{itemize}

\item[{\bf A2}] $\sigma_1^2 = O(n^{\kappa_0})$ for some $\kappa_0 \in
(0,\frac{1}{2})$.

\item[{\bf A3}] $|{\cal A}| = q_n$, $|{\cal B}| = \overline{q}_n$
such that $\overline{q}_n = O(n^{\kappa_1})$ for some $\kappa_1 \in
(0,\frac{1}{2})$.

\item[{\bf A3'}] $p_n$, the number of variables, satisfies the condition that
there is an $\alpha > 0$ such that $\log p_n = O(n^{\alpha})$ for
some $\alpha \in (0,1)$.

\item[{\bf A4}] There exists a $\rho_n$ satisfying
$\rho_n n^{1/2}(\log p_n)^{-1/2} \to \infty$ as $n\to \infty$ such
that
\begin{equation}\label{eq:min_corr}
\min_{j \in {\cal B}}
|\frac{\Sigma_{jy}}{\sqrt{\Sigma_{jj}\sigma_{yy}}}| \geq \rho_n.
\end{equation}

\item[{\bf A5}] There exists a $\overline{\delta}_n$ with $\overline{\delta}_n =
o(\frac{\overline{q}_n}{n \log n})$ such that $\sum_{j\not\in {\cal
B}} \parallel \mathbf{w}_j \parallel_2^2 \leq \overline{\delta}_n$.

\item[{\bf A6}] There exists an $\eta_n > 0$ satisfying
$\eta_n^{-1} = O(n^{\kappa_2})$ for some $\kappa_2 < \frac{1}{2}( 1-
\kappa_0 \vee \kappa_1)$, such that
\begin{equation}\label{eq:min_theta}
\min_{j \in {\cal A}}|\theta_j| \geq \eta_n.
\end{equation}

\item[{\bf A7}] There exists a $\delta \in (0,1)$ such that
\begin{equation}\label{eq:covar_cond}
\parallel\Sigma_{{\cal A}^c {\cal A}} \Sigma_{{\cal A}{\cal A}}^{-1}
\mbox{sign}(\theta_{{\cal A}})\parallel_\infty < \delta.
\end{equation}

\item[{\bf A8}] There is a $\vartheta < \infty$ such that,
\begin{equation}\label{eq:covar_cond_2}
\max_{j \in {\cal A}} \parallel \Sigma_{{\cal A}_j{\cal
A}_j}^{-1}\Sigma_{{\cal A}_j j}
\parallel_1 < \vartheta, ~~\mbox{where}~~{\cal A}_j := {\cal
A}\setminus\{j\}.
\end{equation}
\end{itemize}

A few remarks about these conditions are in order. First, condition
{\bf A1} about the separation of the eigenvalues is not really
necessary, but is assumed to avoid the issue of un-identifiability
of an eigenvector. However, the scaling of the eigenvalues is
important for the analysis. We remark that it is not necessary that
the eigenvalues $\lambda_1,\ldots,\lambda_M$ are the $M$ largest
eigenvalues of $\Sigma$ in order for the conclusions to hold. All
that is necessary is that these are the leading eigenvalues of the
matrix $\Sigma_{{\cal D}{\cal D}}$, and there is enough separation
from the other eigenvalues of $\Sigma$. However, this assumption is
made to simplify the exposition.

Next, the condition that $\overline{q}_n = o(n)$ (implicit from
condition {\bf A3}) is necessary for the consistency of the
estimated eigenvectors $\widetilde{\mathbf{u}}_k$ from Supervised
PCA. Condition {\bf A4} is necessary for the identifiability of the
set ${\cal B}$. {\bf A5} implies that the contribution of the
predictors $\{X_j : j \in {\cal D}\setminus{\cal B}\}$ is negligible
in our analysis. Note that $\overline{\delta}_n$ is essentially
measuring the ``selection bias'' for restricting analysis to ${\cal
B}$ rather than ${\cal D}$. Again, the assumption about the rate of
decay of $\overline{\delta}_n$ can be relaxed at the cost of more
involved analysis and smaller range of values for $\mu_n$ (see also
the remark following \textit{Corollary 1}). Too large a value of
$\overline{\delta}_n$ may mean that we may not be able to select the
variables consistently. Condition {\bf A6} is an identifiability
condition for set ${\cal A}$.



Condition {\bf A7} is needed to guarantee consistency of the
variable selection by LASSO after projection. This condition was
shown to be necessary for variable selection in finite dimensional
LASSO regression by \citeasnoun{zou2005} and also, implicitly by
\citeasnoun{MB2006}. \citeasnoun{ZY2006} termed this the
``irrepresentable condition" and showed that it is nearly necessary
and sufficient for consistency of model selection by LASSO when $p,n
\to \infty$. A sufficient condition for this to hold is that
$\max_{j\in {\cal A}^c}\parallel \Sigma_{{\cal A}{\cal A}}^{-1}
\Sigma_{{\cal A} j}\parallel_1 < \delta$.  Observe that
$\Sigma_{{\cal A}{\cal A}}^{-1} \Sigma_{{\cal A} j}$ is the
population regression coefficient in the regression of $X_j$ on
$\{X_l : l \in {\cal A}\}$. If we are using the estimate $\widehat
\theta^{\widehat{\cal B},\mu}$ then (see proof of \textit{Lemma 2})
we can replace {\bf A7} by the weaker requirement
$$
\parallel\Sigma_{{\cal A}^c \cap {\cal B},{\cal A}} \Sigma_{{\cal A}{\cal
A}}^{-1} \mbox{sign}(\theta_{{\cal A}})\parallel_\infty < \delta,
~~\mbox{for some}~\delta ~\in~ (0,1).
$$

\subsection{LASSO solution}

We use the symbol $\mu$ to denote the penalty parameter in LASSO.
The LASSO estimate of $\theta$, after preconditioning, is given by
\begin{equation}\label{eq:theta_LASSO}
\widehat \theta^\mu = \arg\min_{\zeta \in \mathbb{R}^p} \frac{1}{n}
\parallel \widetilde Y - \mathbf{X}\zeta \parallel_2^2 + \mu
\parallel \zeta\parallel_1.
\end{equation}
We also define the \textit{selected LASSO estimate} of $\theta$ by
\begin{equation}\label{eq:theta_LASSO_B}
\widehat \theta^{\widehat {\cal B},\mu} = \arg\min_{\zeta \in
\mathbb{R}^p, \zeta_{\widehat {\cal B}^c} = 0} \frac{1}{n}
\parallel \widetilde Y - \mathbf{X}\zeta \parallel_2^2 + \mu
\parallel \zeta\parallel_1.
\end{equation}
For future use, we define the \textit{restricted LASSO estimate} of
$\theta$ to be
\begin{equation}\label{eq:theta_LASSO_restr}
\widehat \theta^{{\cal A},\mu} = \arg\min_{\zeta \in \mathbb{R}^p,
\zeta_{{\cal A}^c} = 0} \frac{1}{n}
\parallel \widetilde Y - \mathbf{X}\zeta \parallel_2^2 + \mu
\parallel \zeta\parallel_1.
\end{equation}
The notations used here follow \citeasnoun{MB2006}.

\subsection{Consistency of variable selection}

We shall prove most of our consistency results for the estimate
$\widehat \theta^{\widehat {\cal B},\mu}$ and indicate how (and
under what conditions) the same may be proved for the unrestricted
estimator $\widehat \theta^\mu$. As we shall see, when the model
assumptions hold the former estimator is more reliable under a wider
range of possible dimensions. The latter can consistently select the
model essentially when $p_n = O(n^{\kappa})$ for some $\kappa <
\infty$. In order to prove these results, it will be convenient for
us to assume that we have two independent subsamples of size $n$
each, so that the total sample size is $2n$. And we also assume that
{\bf Step 1} of the variable selection algorithm (estimating ${\cal
B}$) is performed on the first subsample and the other steps are
performed on the second subsample. This extra assumption simplifies
our proofs (see the proof of \textit{Proposition 4} in the Appendix)
somewhat. Further, we shall assume that $K$, the number of latent
components for response $Y$, is known. The results presented here
hold uniformly w.r.t. the parameters satisfying assumptions {\bf
A1}-{\bf A8}.

Let $\widehat {\cal A}_{\widehat{\cal B}, \mu}$ (resp. $\widehat
{\cal A}_\mu$) denote the set of nonzero coordinates of the vector
$\widehat \theta^{\widehat {\cal B},\mu}$ (resp. $\widehat
\theta^\mu$). Whenever the context is clear, we shall drop the
subscripts from $\widehat {\cal A}$. In the following $\zeta$ will
be used to denote a generic value of the parameter.

\vskip.15in \noindent {\bf Proposition 1 :} Let $\widehat {\cal B}$
denote the set of coordinates selected by the preliminary
thresholding scheme of SPC with threshold $\tau_n$. Given any $c_1
> 1$, and there is a $\tau_n(c_1) :=
d_1 \sqrt{\frac{\log p_n}{n}}$, for some constant $d_1 > 2$, such
that, for $n \geq n_{c_1}$,
\begin{equation}\label{eq:SPC_inclusion}
\mathbb{P}(\widehat{\cal B} = {\cal B}) \geq 1 - n^{-c_1}.
\end{equation}

\vskip.1in \textit{Proposition 1} tells us that we can restrict our
analysis to the set ${\cal B}$ while analyzing the effect of
preconditioning, and studying the estimator
$\widehat\theta^{\widehat{\cal B},\mu}$. Our next result is about
the behavior of the estimated eigenvalues and eigenvectors of the
matrix $\mathbf{S}_{\widehat{\cal B}\widehat{\cal B}} :=
\frac{1}{n}\mathbf{X}_{\widehat{\cal B}}^T\mathbf{X}_{\widehat {\cal
B}}$. This result can be proved along the lines of \textit{Theorem
3.2} in \citeasnoun{paul2005}, (see also \citeasnoun{BHPT2006}) and
is omitted.

\vskip.15in \noindent {\bf Proposition 2 :} Let
$(\overline{\mathbf{u}}_{{\cal B}k},\overline{\ell}_k)_{k=1}^K$
denote the first $k$ eigenvector-eigenvalue pairs of $\Sigma_{{\cal
B}{\cal B}}$. Suppose that assumptions {\bf A1}-{\bf A5} hold. Then
there are functions $\gamma_i =
\gamma_i(\lambda_1/\sigma_0,\ldots,\lambda_M/\sigma_0)$, $i=1,2$
such that, given $c_2 > 0$ there exist $d_2,d_2' \geq 1$ so that,
\begin{eqnarray*}
\mathbb{P}(\max_{1\leq k \leq K}
\parallel \widetilde{\mathbf{u}}_{{\cal B}k} -
\overline{\mathbf{u}}_{{\cal B}k}\parallel_2 > d_2 \sigma_0 \gamma_1
\sqrt{\frac{\overline{q}_n \vee \log
n}{n}}(1+\sqrt{\frac{\overline{q}_n \log n}{n}}), ~\widehat{\cal B}
= {\cal B}) &=&
O(n^{-c}),\\
\mathbb{P}(\max_{1\leq k \leq K} |\widetilde \ell_k -
\overline{\ell}_k|
> d_2' \sigma_0^2 \gamma_2 (\sqrt{\frac{\log
n}{n}} + \frac{\overline{q}_n \log n}{n}), ~\widehat{\cal B} = {\cal
B}) &=& O(n^{-c}).
\end{eqnarray*}

\vskip.15in \noindent {\bf Theorem 1 :} Suppose that assumptions
{\bf A1}-{\bf A8} hold. If $\mu = \mu_n$ satisfies $\mu_n =
o(n^{-\kappa_2})$ and $\mu_nn^{\frac{1}{2}(1-\kappa_0\vee \kappa_1)}
\to \infty$ as $n \to \infty$, then there exists some $c
> 1$ such that, for large enough $n$,
\begin{equation}\label{eq:inclusion}
\mathbb{P}(\widehat {\cal A} \subset {\cal A} ) \geq 1 - O(n^{-c}),
\end{equation}
where $\widehat {\cal A} = \widehat {\cal A}_{\widehat{\cal
B},\mu_n}$. If moreover, $p_n$ is such that $\frac{q_n \log p_n}{n}
= o(1)$ as $n \to \infty$, then (\ref{eq:inclusion}) holds with
$\widehat {\cal A} = \widehat {\cal A}_{\mu_n}$.

\vskip.1in \noindent {\bf Theorem 2 :} With $\mu = \mu_n$ and
$\widehat {\cal A}$ as in \textit{Theorem 1}, there exists $c > 1$
such that,
\begin{equation}\label{eq:exclusion}
\mathbb{P}({\cal A}  \subset \widehat {\cal A}) \geq 1 - O(n^{-c}).
\end{equation}

\vskip.1in Clearly, \textit{Theorem 1} and \textit{Theorem 2}
together imply that the SPC/LASSO procedure asymptotically selects
the correct set of predictors under the stated assumptions. The
proofs of these critically rely on the following three results.

\vskip.15in\noindent {\bf Lemma 1 :} Given $\theta \in
\mathbb{R}^p$, let $G(\theta)$ be the vectors whose components are
defined by
\begin{equation}\label{eq:G_theta}
G_j(\theta) = - \frac{2}{n} \langle \widetilde Y - \mathbf{X}\theta,
X_j \rangle
\end{equation}
A vector $\widehat \theta$ with $\widehat \theta_j = 0$ for all $j
\in{\cal A}^c$ is a solution of (\ref{eq:theta_LASSO_restr}) if and
only if, for all $j \in {\cal A}$,
\begin{eqnarray}\label{eq:G_j_theta_cond}
G_j(\widehat \theta) &=&  - ~\mbox{sign}(\widehat \theta_j) \mu ~~
\mbox{if}~\widehat \theta_j \neq 0 \nonumber\\
|G_j(\widehat \theta)| &\leq& \mu ~~\mbox{if}~ \widehat \theta_j = 0
\end{eqnarray}
Moreover, if the solution is not unique and $|G_j(\widehat \theta)|
< \mu$ for some solution $\widehat \theta$, then $\widehat \theta_j
= 0$ for all solutions of (\ref{eq:theta_LASSO_restr}).

\vskip.15in \noindent {\bf Proposition 3 :} Let $\widehat
\theta^{{\cal A},\mu}$ be defined as in
(\ref{eq:theta_LASSO_restr}). Then, under the assumptions of
\textit{Theorem 1}, for any constant $c_3
> 1$, for large enough $n$,
\begin{equation}\label{eq:theta_LASSO_restr_sign}
\mathbb{P}(~\mbox{sign} (\widehat \theta_j^{{\cal A},\mu_n}) =
~\mbox{sign}(\theta_j), ~~\mbox{for all}~j \in {\cal A}) \geq 1-
O(n^{-c_3}).
\end{equation}

\vskip.1in \noindent {\bf Lemma 2 :} Define
\begin{equation}\label{eq:E_mu_B}
{\cal E}_{{\cal B},\mu} = \{ \max_{j \in {\cal A}^c \cap {\cal B}}
|G_j(\widehat \theta^{{\cal A},\mu})| < \mu \} \cap\{ \widehat{\cal
B} = {\cal B}\}
\end{equation}
On ${\cal E}_{{\cal B},\mu}$, $\widehat \theta^{{\cal B},\mu}$ is
the unique solution of (\ref{eq:theta_LASSO_B}) and $\widehat
\theta^{{\cal A},\mu}$ is the unique solution of
(\ref{eq:theta_LASSO_restr}), and $\widehat \theta^{\widehat{\cal
B},\mu} = \widehat \theta^{{\cal A},\mu}$. Also, under the
assumptions of \textit{Theorem 1}, there exists a $c_4
> 1$ such that, for large enough $n$,
\begin{equation}\label{eq:E_mu_bound}
\mathbb{P}({\cal E}_{{\cal B},\mu}^c) = O(n^{-c_4}).
\end{equation}
Further, if we define
\begin{equation}\label{eq:E_mu}
{\cal E}_\mu = \{ \max_{j \in {\cal A}^c} |G_j(\widehat
\theta^{{\cal A},\mu})| < \mu \} \cap\{ \widehat{\cal B} = {\cal
B}\},
\end{equation}
then under the extra assumption that $\frac{q_n\log p_n}{n} = o(1)$,
(\ref{eq:E_mu_bound}) holds with ${\cal E}_{{\cal B},\mu}$ replaced
by ${\cal E}_\mu$. On ${\cal E}_\mu$, $\widehat \theta^{\mu}$ is the
unique solution of (\ref{eq:theta_LASSO}) and $\widehat \theta^{\mu}
= \widehat \theta^{\widehat{\cal B},\mu} = \widehat \theta^{{\cal
A},\mu}$.

\subsection{Effect of projection}

An important consequence of the projection is that the measurement
noise $Z$ is projected onto a $K$ dimensional space (that under our
assumptions also contains the important components of the predictors
of $Y$). This results in a stable behavior of the residual of the
projected response $\Delta$ given by
\begin{equation}\label{eq:Delta}
\Delta := \widetilde Y - \mathbf{X}\theta = \widetilde Y -
\mathbf{X}_{\cal A}\theta_{\cal A}.
\end{equation}
even as dimension $p_n$ becomes large. This can be stated formally
in the following proposition.

\vskip.15in\noindent {\bf Proposition 4 :} Suppose that assumptions
{\bf A1}-{\bf A5} hold. Then there is a constant $\gamma_3 :=
\gamma_3(\sigma_0,\lambda_1,\ldots,\lambda_K+1)$, such that for any
$c_6 > 1$ there exists a constant $d_6 > 0$ so that, for large
enough $n$,
\begin{equation}\label{eq:Delta_bound}
\mathbb{P}(\parallel \Delta \parallel_2 \leq d_6(\gamma_3
\sqrt{\overline{q}_n \vee \log n} + \sigma_1 \sqrt{K \log n})) \geq
1 - n^{-c_6}.
\end{equation}

\vskip.1in As a direct corollary to this we have the following
result about the risk behavior of the OLS-estimator (under $L^2$
loss) of the preconditioned data after we have selected the
variables by solving the optimization problem
(\ref{eq:theta_LASSO_B}).

\vskip.15in\noindent {\bf Corollary 1 :} Suppose that conditions of
\textit{Theorem 1} hold. Then for any $c_7 \geq 1$, there is $d_7 >
0$ such that
\begin{equation}\label{eq:theta_B_risk_bound}
\mathbb{P}(\parallel \widehat \theta^{\widehat{\cal A}_{\widehat
{\cal B},\mu},OLS} - \theta\parallel_2 \leq
d_7\sigma_0^{-1}(\gamma_3 \sqrt{\frac{\overline{q}_n \vee \log
n}{n}} + \sigma_1 \sqrt{\frac{K \log n}{n}})) \geq 1 - n^{-c_7},
\end{equation}
where $\widehat\theta^{\widehat{\cal A}_{\widehat {\cal B},\mu},OLS}
= (\mathbf{X}_{\widehat{\cal A}}^T\mathbf{X}_{\widehat{\cal
A}})^{-1}\mathbf{X}_{\widehat{\cal A}}^T \widetilde Y$, and
$\widehat {\cal A} = \widehat{\cal A}_{\widehat {\cal B},\mu_n} =
\{j \in {\cal P} : \widehat \theta^{\widehat{\cal B},\mu_n} \neq 0
\}$.

\vskip.15in As a comparison we can think of the situation when
${\cal A}$ is actually known, and consider the $L^2$ risk behavior
of the OLS estimator restricted only to the subset of variables
${\cal A}$. Then $\widehat\theta^{{\cal A},OLS} = (\mathbf{X}_{\cal
A}^T\mathbf{X}_{\cal A})^{-1}\mathbf{X}_{\cal A}^T Y$.  Using the
fact that conditional on $\mathbf{X}_{\cal A}$,
$\widehat\theta_{\cal A}^{{\cal A},OLS}$ has $N(\theta_{\cal
A},\sigma_\varepsilon^2 (\mathbf{X}_{\cal A}^T\mathbf{X}_{\cal
A})^{-1})$ distribution, and the fact that the smallest eigenvalue
of $\Sigma_{{\cal A}{\cal A}}^{-1}$ is at least $\ell_1^{-1}$, it
follows (using \textit{Lemma A.1}) that there is a constant $d_7' >
0$ such that
\begin{equation}\label{eq:OLS_A_risk_bound}
\mathbb{P}(\parallel \widehat\theta^{{\cal A},OLS} - \theta
\parallel_2 \geq d_7' \ell_1^{-1/2} \sigma_\varepsilon
\sqrt{\frac{q_n}{n}}) \geq 1 - n^{-c_7}.
\end{equation}
Comparing (\ref{eq:OLS_A_risk_bound}) with
(\ref{eq:theta_B_risk_bound}), we see that if $q_n \gg \log n$ and
$\sigma_1 \gg \sqrt{\overline{q}_n/q_n}$, the estimator
$\widehat\theta^{\widehat{\cal A}_{\widehat {\cal B},\mu},OLS}$ has
better risk performance than $\widehat\theta^{{\cal A},OLS}$.

As a remark, we point out that the bound in
(\ref{eq:theta_B_risk_bound}) can be improved under specific
circumstances (e.g. when $\overline{\delta}_n$, the ``selection
bias'' term defined in {\bf A5}, is of a smaller order) by carrying
out a second order analysis of the eigenvectors
$\{\widetilde{\mathbf{u}}_k\}_{k=1}^K$ (see Appendix of
\citeasnoun{BHPT2006}). The same holds for the bounds on the partial
correlations $\frac{1}{n}\langle (I - P_{\mathbf{X}_{\cal A}} )X_j ,
\widetilde Y \rangle$, for $j \in {\cal A}^c$, given the ``signal''
variables $\{X_l : l \in {\cal A}\}$, that are needed in the proof
of \textit{Proposition 3} and \textit{Lemma 2}. However, the result
is given here just to emphasize the point that preconditioning
stabilizes the fluctuation in $\widetilde Y - \mathbf{X}\theta$, and
so, partly to keep the exposition brief, we do not present the
somewhat tedious and technical work needed to carry out such an
analysis.

As a further comparison, we consider the contribution of the
measurement noise $Z$ in the maximal empirical partial correlation
$\max_{j \in {\cal A}^c} |\frac{1}{n}\langle (I -
P_{\mathbf{X}_{\cal A}} )X_j , \widetilde Y \rangle|$, given $\{X_l
: l \in {\cal A}\}$. For the pre-conditioned response this
contribution is (with probability at least $1-O(n^{-c})$ for some
$c>1$) of the order $O(\frac{\sigma_1 \sqrt{\log n}}{\sqrt{n}})$,
instead of $O(\frac{\sigma_1\sqrt{\log p_n}}{\sqrt{n}})$ as would be
the case if one uses $Y$ instead of $\widetilde Y$. So, if $\log p_n
\gg \log n$, then the contribution is smaller for the
pre-conditioned response. Formalizing this argument, we derive the
following asymptotic result about the model selection property of
LASSO estimator that clearly indicates that under latter
circumstances SPC + LASSO procedure can outperform conventional
LASSO in terms of variable selection.


\vskip.15in \noindent {\bf Proposition 5 :} Suppose that $\log p_n =
c n^{\alpha}$ for some $\alpha \in (0,1)$ and some $c>0$. Suppose
that ${\cal A} = {\cal A}_+ \cup {\cal A}_-$, with ${\cal A}_+$ and
${\cal A}_-$ disjoint and  ${\cal A}_-$ is nonempty
such that $\parallel \theta_{{\cal A}_-}
\parallel_2 = o(n^{-(1-\alpha)/2})$. Assume that $M = K$, ${\cal B}
= {\cal D}$ (so that for all $j \not\in {\cal B}$, $X_j$ are i.i.d.
$N(0,\sigma_0^2)$), and $\sigma_1$ is fixed . Suppose further that
all the assumptions of \textit{Theorem 1} hold,
and there is a $\delta_+  \in (0,1)$ such that (if ${\cal A}_+$ is
nonempty)
\begin{equation}\label{eq:A_plus_condition_1}
\max_{j \not\in {\cal A}_+} \parallel \Sigma_{{\cal A}_+{\cal
A}_+}^{-1}\Sigma_{{\cal A}_+ j} \parallel_1 < \delta_+.
\end{equation}
Then, given $c_8 \geq 1$, for all $\mu_n \geq 0$, for large enough
$n$,
\begin{equation}\label{eq:LASSO_model_selection_failure}
\mathbb{P}(\widehat {\cal A}_{\mu_n}^{LASSO} \neq {\cal A}) \geq 1 -
n^{-c_8},
\end{equation}
where $\widehat {\cal A}_{\mu_n}^{LASSO} = \{j \in {\cal P} :
\widehat \theta_j^{LASSO,\mu_n} \neq 0\}$, where
\begin{equation}\label{eq:LASSO_full}
\widehat \theta^{LASSO,\mu_n} = \arg\min_{\zeta \in \mathbb{R}^p}
\frac{1}{n} \parallel Y - \mathbf{X}\zeta \parallel_2^2 + \mu_n
\parallel \zeta \parallel_1.
\end{equation}


\vskip.15in  {Proposition 5} shows that if $\alpha > 1 -
2\kappa_2$, so that $\eta_n = o(n^{-(1-\alpha)/2})$, and the
assumptions of {Proposition 5} are satisfied, then the SPC +
LASSO approach (solving the optimization problem
(\ref{eq:theta_LASSO_B}) or (\ref{eq:theta_LASSO})) can identify
${\cal A}$ with appropriate choice of penalization parameter $\mu_n$
(as indicated in \textit{Theorem 1}) while LASSO cannot, with any
choice of the penalty parameter.

\section{Classification problems and further topics}
\label{sec.other}

The pre-conditioning idea has potential application in any
supervised learning problem in which the number of features greatly
exceeds the number of observations. A key component is the
availability of a consistent estimator for the construction of the
pre-conditioned outcome variable.

For example, pre-conditioning can be applied to classification
problems. Conceptually, we separate the problems of a) obtaining a
good classifier and b) selecting a small set of good features for
classification. Many classifiers, such as the support vector
machine, are effective at finding a good separator for the classes.
However they are much less effective in distilling these features
down into  a smaller set of uncorrelated features.

Consider a two-class problem, and suppose we have trained a
classifier, yielding estimates $\hat p_i$, the probability of class
2 for observation $i=1,2,\ldots N$. Then in the second stage, we
apply a selection procedure such as forward stepwise or the LASSO,
to an appropriate function of $\hat p_i$; the quantity $\log[\hat
p_i/(1-\hat p_i)]$ is a logical choice.

We generated data as in  example of section \ref{sec.example};
however we turned it into a classification problem by defining the
outcome class $g_i$ as 1 if $y_i<0$ and 2 otherwise. We applied the
nearest shrunken centroid (NSC) classifier of \citeasnoun{THNC2002},
a  method for classifying microarray samples. We applied
forward stepwise regression both to $g_i$ directly (labeled FS), and
to the output $\log(\hat p_i/(1-\hat p_i))$ of the NSC classifier
(labeled NSC/FS).

\begin{figure}
\centerline{\epsfig{file=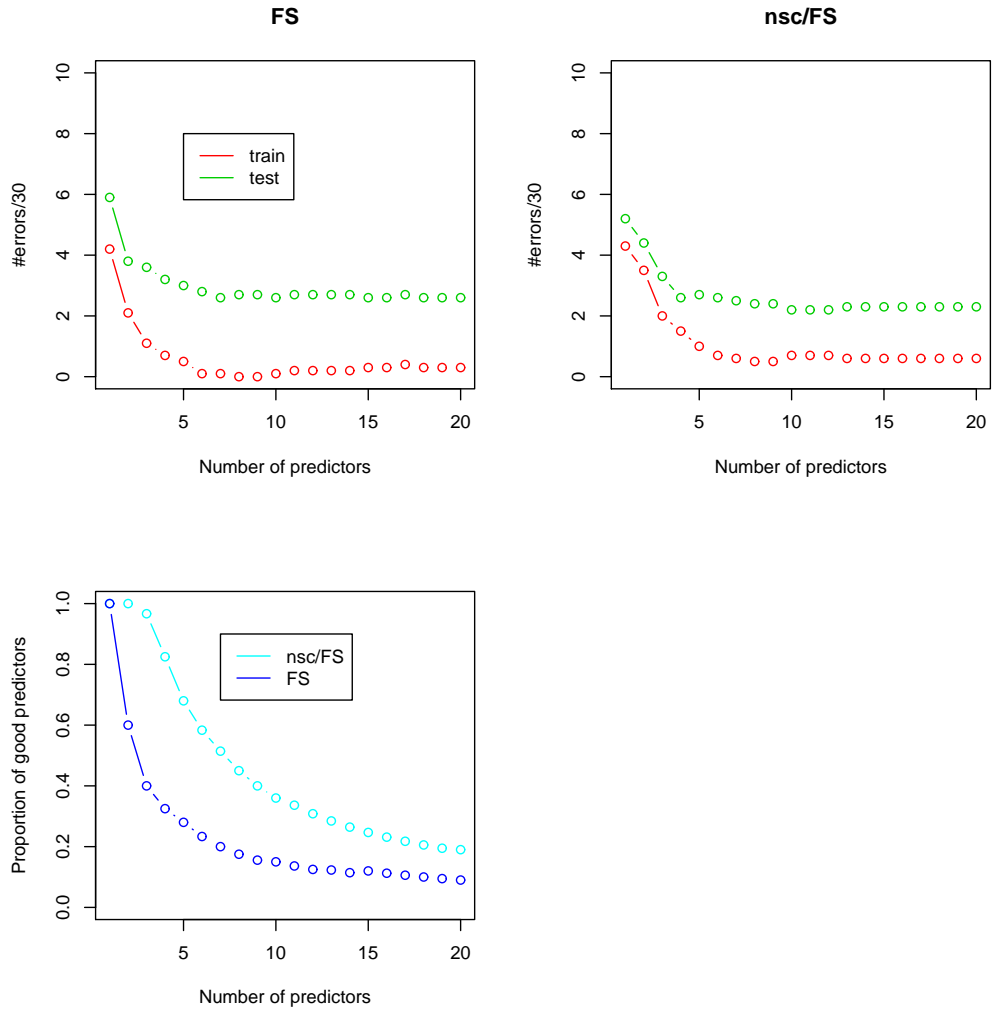,width=\textwidth}}
\caption[figclass]{\small\em Results of applying pre-conditioning in
a classification setting. Top left panel shows teh number of test
misclassification errors from  forward stepwise regression; in teh
top right panel we have applied forward stepwise regression to the
pre-conditioned estimates from nearest shrunken centroid classifier.
The proportion of good predictors selected by each method is shown
in the bottom left.} \label{figclass}
\end{figure}

The results of 10 simulations are shown in Figure \ref{figclass}. We
see that NSC/FS does not improve the test error of FS, but as shown
in the bottom left panel, it does increase the number of ``good''
predictors that are found. This is a topic of further study.

\vskip.2in \noindent {\bf Acknowledgments:}

We thank the referees and editors for comments that led to
improvements in this work.
Hastie was partially supported by grant DMS-0505676 from the
National Science Foundation and grant 2R01 CA 72028-07 from the
National Institute of Health. Tibshirani was partially supported by
National Science Foundation Grant DMS-9971405 and National
Institutes of Health Contract N01-HV-28183.

\section*{Appendix}
A full version of this paper that includes the Appendix is available at
\begin{verbatim}
http://www-stat.stanford.edu/~tibs/ftp/precond.pdf
\end{verbatim}
 and also  in  arXiv  archive.

\bibliography{tibs,/home/hastie/bibtex/tibs.bib}

\bibliographystyle{agsm}

\end{document}